\documentclass{amsart}
\usepackage{amssymb,latexsym,verbatim,hyperref}
\usepackage[all]{xy}

\newcommand\DMO[2]{\DeclareMathOperator{#1}{#2}}
\newcommand\comp{\circ}

\newcommand\Zn{\mathbb{Z}}

\newcommand\Nn{\mathbb{N}}
\newcommand\Rn{\mathbb{R}}
\newcommand\CP{\mathbb{C}\mathrm{P}}

\newcommand\ckH{{\check{H}}}

\newcommand\qbound{\preccurlyeq}
\newcommand\qequiv{\approx}

\newcommand\homot{\simeq}
\newcommand\bdy{\partial}
\renewcommand\bar[1]{\overline{#1}}
\newtheorem{theorem}{Theorem}
\newtheorem{corollary}{Corollary}
\newtheorem{lemma}{Lemma}
\theoremstyle{definition}
\newtheorem{definition}{Definition}

\DMO{\Obj}{Obj} \DMO{\Mor}{Mor} \DMO{\im}{im} \DMO{\coker}{coker}
\DMO{\Hom}{Hom} \DMO{\Spec}{Spec} \DMO{\Frac}{Frac} \DMO{\tr}{tr}
\DMO{\Sym}{Sym} \DMO{\supp}{supp} \DMO{\Proj}{Proj} \DMO{\vol}{vol}
\DMO{\Riem}{Riem} \DMO{\Diff}{Diff} \DMO{\FV}{FV} \DMO{\Vol}{Vol}
\DMO{\dist}{dist} \DMO{\Cont}{Cont}

\newcommand\Volcur{\mathbf{M}}
\newcommand\Volmet{\Vol_\mathrm{met}}
\newcommand\Volch{\Vol_\mathrm{ch}}
\newcommand\VolAWP{\Vol_\mathrm{AWP}}
\newcommand\Voltr{\Vol_\mathrm{tr}}
\newcommand\Volcoh{\Vol_\mathrm{coh}}

\newcommand\FVcur{\FV_\mathrm{\!cur}}
\newcommand\FVmet{\FV_\mathrm{\!met}}
\newcommand\FVch{\FV_\mathrm{\!ch}}
\newcommand\FVcell{\FV_\mathrm{\!cell}}

\newcommand\Phicur{\Phi_\mathrm{cur}}
\newcommand\Phimet{\Phi_\mathrm{met}}
\newcommand\Phich{\Phi_\mathrm{ch}}
\newcommand\Phicell{\Phi_\mathrm{cell}}

\renewcommand\tilde{\widetilde}

\begin{document}

\title{Generalized Dehn functions I}
\author{Chad Groft}
\date{January 15, 2009}

\begin{abstract}
  Let $X$ be a finite CW complex or compact Lipschitz neighborhood
  retract, and let $\tilde X$ be its
  universal cover; let $M$ be a compact orientable manifold of
  dimension $q\ge 2$ and $\bdy M\ne\emptyset$.  We establish the
  existence of an isoperimetric profile for functions $f\colon
  M\to\tilde X$, in the metric and cellular senses, and show that they
  are equivalent up to scaling factors when $X$ is a triangulated CLNR
  (for example a triangulated Riemannian manifold).  We also show that
  two finite complexes $X$ and $Y$ have the same profiles up to
  scaling given the existence of a $q$-connected map between them.
\end{abstract}
\maketitle

\section*{Introduction}\label{S:intro}

Let $G$ be a group with finite presentation $P=\left< X|R\right>$,
so that there is an epimorphism $\pi\colon F(X)\to G$ whose kernel
is the normal subgroup of $F(X)$ generated by $R$. If $w\in F(X)$
is a word such that $\pi(w)=e$, then
$w=(r_1^{w_1})^{\epsilon_1}\dots(r_N^{w_N})^{\epsilon_N}$
for some $N\ge 0$, $r_i\in R$, $w_i\in F(X)$, and $\epsilon_i=\pm1$.
The smallest possible $N$ is called the \emph{filling volume}
for $w$, or $\FV(w)$, and the \emph{Dehn function} for the presentation
$P$ is $\Phi_P\colon\Nn\to\Nn$,
\[
  \Phi_P(n) = \sup\{\,\FV(w) : \pi(w)=e,\,|w|\le n\,\},
\]
where $|w|$ is the word length of $w$ in the generators $X$. $\Phi_P$
is not exactly an invariant of $G$; for example, the trivial group has
presentations $P=\left<|\right>$ and $Q=\left<x|x\right>$, and
$\Phi_P\equiv 0$ while $\Phi_Q(n)=n$ for all $n$. However, for any
two finite presentations $P$ and $Q$ of the same group, we have
\[
  \Phi_P(n) \le A\cdot\Phi_Q(Bn) + Cn + D
\]
for some $A$, $B$, $C$, $D$, and \emph{vice versa}.
This is an equivalence relation between functions, so---in the typical
abuse of language---one often speaks of the Dehn function of $G$.
Dehn functions can vary widely, from linear (as with $G$ trivial)
to busy-beaver (as when $G$ has word problem equivalent to the
halting problem).

Burillo and Taback give a geometric interpretation to the Dehn function in
\cite{jB02}: Fix a compact
Riemannian manifold $M$. Let $f\colon[0,\infty)\to[0,\infty)$ be the
smallest function where, given a loop
$\gamma\colon S^1\to\tilde M$ of length $\ell$, there is a disk map
$\Gamma\colon D^2\to\tilde M$ where $\Gamma\restriction S^1=\gamma$
(we say that $\Gamma$ \emph{fills} $\gamma$)
and where the area of $\Gamma$ is at most $f(\ell)$. Then
$f$ is equivalent to the Dehn function of $\pi_1(M)$. This tells us that part of
the large-scale geometry of $\tilde M$ is determined by $\pi_1(M)$.

The argument is as follows: Triangulate $M$, and extract a presentation
$\left<X|R\right>$ for $G=\pi_1(M)$ from the triangulation. An word
$w\in F(X)$ is represented by a loop $\gamma_w$ in the 1-skeleton of $M$,
whose length is roughly a constant times $|w|$; $\gamma_w$ lifts to $\tilde M$
iff it is nullhomotopic, which occurs iff $\pi(w)=e$. Conversely, a loop
$\gamma$ in $M$ may be deformed to the 1-skeleton, where
it represents a word $w_\gamma\in F(X)$, and $\gamma$ lifts to $\tilde M$
iff $\pi(w_\gamma)=e$.

If a loop $\gamma$ is nullhomotopic, then
$w_\gamma=(r_1^{w_1})^{\epsilon_1}\dots(r_N^{w_N})^{\epsilon_N}$, and
there is a nullhomotopy of $\gamma$ which deforms $\gamma$ to the
1-skeleton of $\tilde M$, then traverses $N$ 2-cells
of $\tilde M$ corresponding to the words $r_i^{w_i}$; the $\epsilon_i$
determine the direction in which the cell is traversed. This fact, combined
with the finite area of the 2-cells, bounds $f$
in terms of $\Phi_G$. Conversely, if $\pi(w)=e$, then $\gamma_w$ is
nullhomotopic, and the nullhomotopy may be deformed to a disk map
$\Gamma\colon D^2\to\tilde M^{(2)}$. Now $\Gamma$ represents an
element $[\Gamma]\in \pi_2(M^{(2)},M^{(1)})$, whose generators as a group
correspond both to 2-cells in $M$ and to words $r^v$ where $r\in R$
 and $v\in F(X)$. Since $\bdy[\Gamma]=[\Gamma\restriction S^1]=[\gamma_w]=w$,
we can bound $\FV(w)$, and therefore $\Phi_G$, in terms of $f$.

It should be noted that $M$ can be replaced by an arbitrary
CW complex with finite 2-skeleton whose cells are given unit volume
by convention.

Alonso, Wang, and Pride extend this idea to higher dimensions
in \cite{jmA99}. Say that a group $G$ is of type $F_q$
if there is a CW complex $X=K(G,1)$ with finite $q$-skeleton.
Thus $F_0$ is trivial, $F_1$ means ``finitely generated'', $F_2$
means ``finitely presented'', and so forth. For such an $X$,
and for $q\ge 2$,
there is a $q$-dimensional Dehn function which compares the
volumes of cellular maps $f\colon S^{q-1}\to X$ to those of cellular
extensions $g\colon D^q\to X$. (Here ``volume'' is interpreted as word
length in, say, $\pi_q(X^{(q)},X^{(q-1)})$; thus the case $q=2$ gives us
the classical Dehn function.) Alonso \emph{et.~al.}\ prove
that this higher-dimensional Dehn function is independent of the
specific $K(G,1)$ chosen, up to equivalence.

In \cite[section 6.D]{mG01} Gromov correctly asserts, but does not prove, the following:
Let $M$ be a compact Riemannian manifold with
$\pi_2(M)=\dots=\pi_{q-1}(M)=0$, so that $M$ is ``approximately'' a $K(G,1)$.
Let $f$ be the isoperimetric function for filling
maps $S^{q-1}\to\tilde M$ with maps $D^q\to\tilde M$.
Then $f$ is equivalent to the $q$-dimensional Dehn function on
$\pi_1(M)$. Again we see an aspect of the large-scale geometry of
$\tilde M$ which depends solely on the low-dimensional topology
of $M$.

Gromov continues in \cite[remark $6.34{\frac12}_+$(c)]{mG01}:
\begin{quotation}
  It is unclear how to approach the $q$-dimensional isoperimetric
  problem with $q\ge 3$ without the \dots condition $\pi_i=0$ for
  $i\le q-1$\dots.  Also, one should be more specific here on the
  allowed topology of manifolds \dots involved in the definition.  For
  example, one can stick to $(q-1)$-spheres which are allowed to be
  filled in (spanned) by $q$-balls.  Or, we can look at $(q-1)$-tori
  filled by solid $q$-tori\dots.
\end{quotation}
In general, we could look at maps defined on any manifold $N^q$ with
nonempty boundary; for convenience, we only consider $N$ connected
and orientable, though a parallel theory for $N$ non-orientable
doubtless exists. Specifically, given a map $f\colon \bdy N\to\tilde M$,
extend to a map $N\to\tilde M$ with as small a volume as possible,
then bound this minimal volume in terms of the volume of $f$.
Such functions are examined in \cite{nB09}, for example.

As a
sort of limiting case, one can also consider $(q-1)$-currents filled by
$q$-currents; this version is seen in \cite{dbaE92} and \cite{mG93} in the case of highly
connected spaces.

The question is in fact nontrivial if we allow $\tilde M$
not highly connected. This is most easily seen by considering $M=\CP^2$ under,
say, the Fubini-Study metric. Topologically, $M$ is a CW complex
consisting of a 4-cell attached to a 2-sphere by the Hopf map
$\eta\colon S^3\to S^2$; this map represents a generator of $\pi_3(S^2)=\Zn$.
If $f\colon S^3\to\CP^2$ is the composition of a representative
map for $n[\eta]$ with the inclusion $S^2\to\CP^2$, then $f$ has zero
volume, but any filling map $D^4\to\CP^2$ must have volume at least
$n\Vol(\CP^2)$. Thus the isoperimetric function corresponding to maps
$S^3\to\CP^2$ and $D^4\to\CP^2$ is infinite everywhere, whereas looking
only at $\pi_1(\CP^2)=e$, one would expect a finite-everywhere, even
linear, function. If one replaces $D^4$ with $\CP^2\setminus B$ where $B$
is a small open ball, one can fill $f$ with a map of volume $(n-1)\Vol(\CP^2)$,
but one still obtains an infinite isoperimetric function. Contrast with the
version involving chains; there are no nonzero cellular 3-chains, so the
cellular isoperimetric
function is zero everywhere. We therefore expect that a 3-dimensional
current $S$ can be filled with a 4-dimensional current $T$ where
$\Volcur(T) \le c\Volcur(S)+d$ with $c$, $d$ constant; and indeed this occurs.

In this paper we establish technical results for isoperimetric profiles for maps
$M^q\to\tilde X$, where $M$ is a fixed compact manifold
with nonempty boundary of dimension $q\ge 2$, and $X$ is a fixed
compact Riemannian manifold.

In particular we show that, up to equivalence, this profile depends only on the
``low-dimensional topology'' of $X$. Specifically, let $Y$ be another
compact Riemannian manifold, and let $f\colon X\to Y$ be a continuous
function where $f_*\colon \pi_t(X)\to\pi_t(Y)$ is an isomorphism
for $1\le t<q$. Then $X$ and $Y$ have the same profiles for any
fixed $M$. This idea is borrowed from \cite{jmA99}, in which such a function
is constructed for $X$ and $Y$ both $K(G,1)$'s for a common group $G$.
Such spaces are homotopy equivalent, as is well known. This
equivalence, along with a bound on the local volume of the homotopies,
is how Alonso \emph{et.\ al.\/}\ establish the invariance of their
Dehn functions. We extend their idea by noting that since only the
$q$-skeleta of $X$ and $Y$ matter, the homotopy groups $\pi_t$
for $t\ge q$ can be changed essentially at will and must be irrelevant.
In particular, this implies that these generalized Dehn functions are
well-defined, up to equivalence, on groups $G$ of type $F_q$.

We also establish homological versions of these profiles, in which we
replace maps from $M$ with $q$-currents or $q$-chains, and prove
similar theorems about them.

As when $M=D^q$ and $\tilde X$ is
$(q-1)$-connected, this is done by addressing a similar question
where $X$ is a CW complex, then connecting the two notions
through a variant of the Deformation Theorem.
There is one serious issue with doing so: given a function
$f\colon M\to X$, how does one determine the volume of $f$?
To have a volume in terms of number of cells covered in some
sense, we must at least have $f[M]\subseteq X^{(q)}$ and
$f[\bdy M]\subseteq X^{(q-1)}$, as in \cite{jmA99} for the cases
$M=D^q$ and $M=S^q$.

In the literature (for example \cite{nB09}), one sees the notion of an
\emph{admissible} map, \emph{i.e.}, a map $f$ as above and for
which $f^{-1}[X\setminus X^{(q-1)}]$ is a disjoint union of open
$q$-dimensional disks, each of which is mapped homeomorphically
onto the interior of a $q$-cell of $X$. The volume of $f$,
hereafter written $\Vol f$, is
simply the number of disks. While many maps are admissible
(certainly enough to define the relevant profiles), not all are.
In particular, if a Lipschitz map $f\colon M\to X$ has been
deformed to a map $f'\colon M\to X^{(q)}$, it is unlikely
that $f'$ will be admissible.

If $f$ happens to be cellular for some triangulation of $M$,
then we may restrict to each $q$-cell $\Delta$ of $M$, take
the word length in $\pi_q(X^{(q)},X^{(q-1)})$ as the volume
of $f\restriction\Delta$ (as is done in \cite{jmA99}),
and sum over all $\Delta$. This is very useful for
explicitly constructed maps, but again it is a bit much to ask from an
arbitrary map.

Alternately, generalize the concept of ``admissible'' somewhat, and
suppose that each component $W$
of $f^{-1}[X^{(q)}\setminus X^{(q-1)}]$, though not necessarily a disk, has smooth
boundary. Then, ignoring some boundary issues, $f\restriction\bar W$
is a map from one orientable $q$-dimensional manifold with boundary
to another (namely a $q$-cell of $X$). This map has a certain degree, taken
as nonnegative; we can then sum this degree over all $W$.
A variant of this idea, which uses \v Cech cohomology,
can be applied to any map $f\colon M\to X^{(q)}$
with $f[\bdy M]\subseteq X^{(q-1)}$.

We spend a fair amount of effort showing that these definitions are
equivalent for our purposes.

The paper is organized as follows.  Section~\ref{S:definitions}
defines the notions of volumes for maps and currents in our various
spaces. Section~\ref{S:hinv} is concerned with insuring
that the cellular definitions of volume for a map are equivalent
up to homotopy, which allows us to formally define
the filling volumes and isoperimetric profiles in section~\ref{S:filling}.
We relate the metric and cellular profiles and notions of
volume in section~\ref{S:metvscell}; if $X$ is a compact
Riemannian manifold with a triangulation, then
the profiles for a covering space of $X$ are the same up to equivalence.

Finally, we show in section~\ref{S:homotopy} the cellular version
of the result above:
that if a map $f\colon X\to Y$ exists which is a $q$-homotopy
equivalence, then the $q$-dimensional profiles for $X$ and $Y$
are equivalent.  In fact, the result
is slightly stronger; the induced map $f_*\colon \pi_q(X)\to\pi_q(Y)$
need not be onto.  This may be seen as generalizing the result given
in \cite{jmA99}, where the $q$-dimensional profiles on a highly
connected $\tilde X$ are shown to depend essentially on $\pi_1(X)$
alone.

If $X$ is a compact space with $\pi_1(X)=G$ and where $\tilde X$ is
$(q-1)$-connected, one can apply these constructions and theorems to
establish geometric group invariants for $G$ which are apparently
distinct from the $q$-dimensional Dehn function.  In part~II of the
paper, we will
see that the new profiles are equal for $X$ a CW complex, and
equal almost everywhere for $X$ a manifold, provided $q\ge
4$.  In the case $q=3$ the profile obtained is dependent only on $\bdy
M$.  For $q=2$ we have no positive results.

Except for some corrections and technical additions, this work has been published
as a dissertation (see \cite{clG07}). The author would like to thank
Richard Schoen, under whose advice the dissertation was completed,
and Alex Nabutovsky and Rina Rotman, for valuable criticism of the articles.

\section{Definitions and technical lemmas}\label{S:definitions}

\subsection{Conventions}

For the cellular definitions and results, $X$ will denote a
connected CW complex. We assume in that each attaching
map is cellular, or equivalently that there
is a vertex in the range of each attaching map.
For the metric definitions and results, $X$ will denote a
connected local Lipschitz neighborhood retract, or LLNR.
(That is, $X$ will be a subset of some
$\Rn^N$ with a neighborhood $U$ and a locally Lipschitz retraction
$r\colon U\to X$. This is slightly more general than a manifold,
and is the natural setting for most of our metric work.)
In many cases $X$ will be an LLNR with a Lipschitz
triangulation, putting it in both categories.

$q$ will always denote a positive natural number, and $M$ will always
be a compact orientable smooth manifold of dimension $q$, possibly
with boundary.  If $\bdy M$ is nonempty, we will assume $M$ is
connected.  In the case $q=1$, $\bdy M$ will always be empty, so that
$M$ is a finite disjoint union of $S^1$s.  We assume that $M$ has a
basepoint $*$, and that $\bdy M$ contains $*$ unless it is empty.

For convenience, we define $\bdy f := f\restriction\bdy M$ for any
$f\colon M\to X$.

\subsection{Volumes}

When $X$ is an LLNR, the volumes of currents in $X$ and Lipschitz maps
to $X$ are well established.  If $T$ is an integral $q$-current in
$X$, we take its mass $\Volcur(T)$ as defined in \cite[\S 4.1.7]{hF69}
as the relevant volume.  If $f\colon M\to X$ is a Lipschitz map under
some (any) metric on $M$, we define its volume as usual:
\begin{equation*}
  \Volmet f :=
  \int_M \sqrt{\det(f^*g)_{ij}}\ dx^1\wedge\dots\wedge dx^q =
  \int_M |J\!f(x)|\,dV_x
\end{equation*}
where $(x^1,\dots,x^q)$ are local coordinates.

Now let $X$ be a CW complex.  Here, all our volumes will be word
lengths in various groups.  Recall that, if $G$ is a group with
generating set $S$, the \emph{length} of $g\in G$ (which we denote
$\|g\|$) is the least $n$ where $g=s_1^{\epsilon_1}\dots
s_n^{\epsilon_n}$ for some $s_i\in S$ and some $\epsilon_i=\pm1$.  By
obvious convention, $\|e\|=0$.

Let $C_q(X)=C_q^{CW}(X) = H_q(X^{(q)},X^{(q-1)})$ be the free abelian
group with basis $S$ in 1-1 correspondence with the $q$-cells of $X$.
The map $\bdy\colon C_q(X)\to C_{q-1}(X)$ is canonical.  We define the
volume of $c\in C_q(X)$ as the word length $\|c\|$ in the generators
$S$.

As noted above, we only assign volumes to certain functions $f\colon M\to X$.

\begin{definition}
  A continuous function $f\colon M\to X$ is \emph{quasi-cellular} if
  it is a map $f\colon (M,\bdy M,*)\to(X^{(q)},X^{(q-1)},X^{(0)})$.
  For such a function $f$, we define the \emph{content} of $f$ as
  the set $\Cont(f) := M\setminus f^{-1}[X^{(q-1)}]$.
\end{definition}

If $f$ is quasi-cellular, so is $\bdy f$.  The basepoint condition is
for convenience.  The intuition behind $\Cont(f)$ is that the part of
$f$ mapping to the $(q-1)$-skeleton cannot contribute to the volume of
the $q$-dimensional object $f$.

One might consider taking the volume of a quasi-cellular map $f$ to be
the length of $f_*([M])\in C_q(X)$, but this reduces to the volume of
a chain, which we are trying to avoid.  If it is necessary to cover
the same cell in opposite orientations in order to fill a given map on
$\bdy M$, we want both of these to have positive volume, not cancel
each other out.

\subsubsection{Alonso-Wang-Pride volume}
In the case $M=D^q$ or $M=S^q$ (which by the collapsing map
$(D^q,S^{q-1},*)\to(S^q,*,*)$ may be seen as a special case of the
first), a definition already exists, introduced by Alonso, Wang, and Pride
in \cite{jmA99}.  They define the groups
\begin{alignat*}{9}
  K_t(X,v) &= \pi_t(X^{(t)},X^{(t-1)},v), \qquad & t\ge 2, \\
  K_1(X) = K_1(X,v) &= \pi_1(X^{(1)}/X^{(0)},*), & t=1,
\end{alignat*}
where $v$ is a 0-cell of $X$.  For $t\ge 3$, $K_t(X,v)$ is a free
$\pi_1(X^{(t-1)},v)$-module, so it is a free abelian group with
basis
\begin{equation*}
  S = \{\, \gamma * \phi_\sigma : \gamma\in\pi_1(X^{(t-1)},v),\,
  \textrm{$\sigma$ a $t$-cell}\,\}
\end{equation*}
where $\phi_\sigma$ is the homotopy class of the standard map
$i\colon D^t\to X^{(t)}$ combined with some fixed
path from a 0-cell on $\bdy\sigma$ to $v$.  The group $K_t(X,v)$
is independent of $v$ up to isomorphisms which preserve the set of
generators and therefore the word length.

$K_2(X,v)$ is not generally abelian, but it is a free crossed
$\bigl(\pi_1(X^{(1)},v),\bdy\bigr)$-module as in \cite{jhcW49b}, and
there is a set of group generators
\begin{equation*}
  S = \{\,\gamma*\phi_\sigma : \gamma\in\pi_1(X^{(1)},v),\,
  \textrm{$\sigma$ a 2-cell}\,\}
\end{equation*}
defined as above.  $X^{(1)}$ is path-connected, and associated to each
path $\gamma\colon v\rightsquigarrow w$ is an isomorphism
$p_\gamma\colon K_2(X,v)\to K_2(X,w)$ which is a bijection between the
sets of group generators.  This is enough to make the definition of
volume independent of $v$, as noted in Lemma~1 of \cite{jmA99}.

$K_1(X)$ is a free group with generators corresponding to the 1-cells
of $X$.

\begin{definition}
  If $f\colon D^q\to X$ or $f\colon S^q\to X$ is quasi-cellular, we
  define the \emph{AWP-volume of $f$} as $\VolAWP f:=\|[f]\|$, where
  $[f]$ is the equivalence class of $f$ in $K_q(X,f(*))$.  If $M=S^1$,
  then take $\VolAWP f := \|[\pi\comp f]\|$ where
  $\pi$ is the natural projection
  $(X^{(1)},X^{(0)})\to(X^{(1)}/X^{(0)},*)$.
\end{definition}

\subsubsection{Triangulated volume}
Return to general $M$, and let $\tau$ be a $C^1$ triangulation
of $M$ in the sense of \cite{jrM63}. Assume $*\in\tau^{(0)}$.

\begin{definition}
  A continuous map $f\colon M\to X$ is \emph{$\tau$-cellular} if it is
  a map
  \begin{equation*}
    f\colon(M,\tau^{(q-1)},*)\to(X^{(q)},X^{(q-1)},X^{(0)}).\footnote{This
      should really be called``$\tau$-quasi-cellular'', but that would be
      ugly.}
  \end{equation*}
\end{definition}

A $\tau$-cellular map is automatically quasi-cellular.

If we restrict a $\tau$-cellular map $f$ to a single $q$-simplex of
$\tau$, we obtain a quasi-cellular map $D^q\to X$, which has an
AWP-volume.  By summing over the top-dimensional simplices of $\tau$,
we obtain a reasonable definition of volume for $f$.

\begin{definition}
  For $\tau$ a $C^1$ triangulation of $M$ and $f\colon M\to X$ a
  $\tau$-cellular map, the \emph{$\tau$-volume of $f$} is
  \begin{equation*}
    \Vol_\tau f := \sum_\Delta \VolAWP(f\restriction\Delta),
  \end{equation*}
  where $\Delta$ ranges over the top-dimensional simplices of $\tau$.
\end{definition}

\begin{definition}
  For $f\colon M\to X$ quasi-cellular, the \emph{triangulated volume}
  of $f$ is
  \begin{equation*}
    \Voltr f := \inf_\tau \Vol_\tau f,
  \end{equation*}
  where $\tau$ ranges over those triangulations such that $f$ is
  $\tau$-cellular.
\end{definition}

As we will see, we will often construct $f$ cell by cell on a fixed
triangulation.  For such $f$, this notion of volume is convenient.
However, $\Voltr f$ is finite iff there is a triangulation $\tau$
where $\Cont(f)\cap\tau^{(q-1)}=\emptyset$, which fails for a number
of maps.  Infinity is a poor reflection of $\Volmet f$ in such cases.

\subsubsection{Cohomology volume}
Finally, we turn to (co)homology.  Recall that
our objection to using $\|f_*([M])\|$ as the volume of $f$ was that
different regions of $M$ might overlap the same point with opposite
orientations.  We can minimize this issue by considering each
component $W$ of $\Cont(f)$ separately.  Assuming that $H_q(\bar
W,\bdy W)\cong\Zn$ for each $W$, we can consider the degree of
\begin{equation*}
f_*\colon H_q(\bar W,\bdy W)\to H_q(X^{(q)},X^{(q)}\setminus\sigma(W))
\end{equation*}
where $\sigma(W)$ is the
unique $q$-cell of $X$ which meets $f[W]$.  The volume could then be a
sum of degrees over components of $\Cont(f)$.

This scheme has two potential problems.  One is that there may be
cancellation within a single component $W$.  If $W$ and $f\restriction W$
are sufficiently nice, we can deal with this, as we will see in the proof of lemma
\ref{L:trvscoh}.  A more fundamental problem is that we
generally do not know that $H_q(\bar W,\bdy W)$ is infinite cyclic.
To fix this problem, we replace homology by \v Cech cohomology, for
which \cite[Ch. XI, \S6]{sE52} is the standard reference.

If $\sigma$ is a $q$-cell of $X$, then
$\ckH^q(X^{(q)},X^{(q)}\setminus\sigma^\circ)\cong\Zn$.  If
$W\subseteq M^\circ$ is open and $\{\,W_\alpha : \alpha\in I\,\}$ is
an indexed collection of its components, then each
$\ckH^q(M,M\setminus W_\alpha)$ is infinite cyclic and the diagram
$\ckH^q(M,M\setminus W_\alpha)\to \ckH^q(M,M\setminus W)$,
 where $\alpha$ ranges over $I$, is a
coproduct diagram.  It follows that
$\ckH^q(M,M\setminus W)$ is a free abelian group with basis $S\cong
I$; the coordinate corresponding to $\alpha$ can be seen as the number
of times that $W_\alpha$ covers $\sigma$.

\begin{definition}\label{D:volcoh}
  For a quasi-cellular map $f\colon M\to X$, and for each $q$-cell
  $\sigma$ in $X$, there is a map
  \begin{equation*}
    f^*_{(\sigma)}\colon \Zn \cong \ckH^q(X^{(q)},X^{(q)}\setminus\sigma^\circ)
    \to \ckH^q(M,M\setminus f^{-1}[\sigma^\circ]) \cong
    \bigoplus_{\alpha\in I_\sigma} \Zn,
  \end{equation*}
  where $I_\sigma$ indexes the components of $f^{-1}[\sigma^\circ]$.
  The \emph{cohomology volume} of $f$ is
  \begin{equation*}
    \Volcoh f := \sum_\sigma \| f^*_{(\sigma)}([X^{(q)}]) \|,
  \end{equation*}
  where $\sigma$ ranges over the $q$-cells of $X$ and $[X^{(q)}]$
  generates $\ckH^q(X^{(q)},X^{(q)}\setminus\sigma^\circ)$.
\end{definition}

Note that $\Volcoh f$ is always finite, since $f[M]$ is compact and
therefore meets the interior of only finitely many $q$-cells (see
\cite[theorem A.1]{aH02}). Also note that, if $f$ is admissible,
then $\Vol f=\Volcoh f$.

\section{Homotopy ``invariance'' of volume}\label{S:hinv}

Let $X$ be a CW complex for this section.  It is annoying to have two
or three notions of volume for a single quasi-cellular map $f\colon
M\to X$.  Of course one cannot say in general that $\VolAWP(f) =
\Voltr(f) = \Volcoh(f)$.  However, we have a substitute.

\begin{definition}\label{D:nicehomotopy}
  Let $f,\,f'\colon M\to X$ be quasi-cellular.  A homotopy $H\colon
  f\homot f'$ is \emph{nice} if it fixes $\bdy M$ and $*$ and if
  $H[M\times[0,1]]\subseteq X^{(q)}$.  $f$ and $f'$ are called
  \emph{nicely homotopic} if such an $H$ exists.
\end{definition}

Note that for $f,\,f'\colon M\to X$ where $M=D^q$ or $M=S^q$, if $f$
and $f'$ are nicely homotopic, then $\VolAWP f=\VolAWP f'$.

This section is devoted to these two lemmas:

\begin{lemma}\label{L:AWPvscoh}
  Let $M=D^q$ or $M=S^q$, and let $f\colon M\to X$ be quasi-cellular.
  Then $\VolAWP f \le \Volcoh f$, and $\Volcoh g\le \VolAWP f$
  for some $g$ which is nicely homotopic to $f$ and admissible.
\end{lemma}

\begin{lemma}\label{L:trvscoh}
  Let $f\colon M\to X$ be quasi-cellular.  Then $f$ is nicely
  homotopic to $g$ where $\Voltr g = \Vol_\tau g \le \Volcoh f$ for some
  $\tau$ where, for each $q$-cell $\Delta$, $g[\Delta]$ meets the interior
  of at most one cell $\sigma$; if there is such a $\sigma$, then $g$ is
  smooth on $g^{-1}[\sigma^\circ]$ and covers each point of
  $\sigma^\circ$ at most $|\deg(g\restriction\Delta)|$ times. Also $f$ is nicely
  homotopic to $h$ where $\Volcoh h \le \Voltr f$ and $h$ is admissible.
\end{lemma}

Note that any quasi-cellular map $f\colon M\to X$ actually maps into a
finite subcomplex $Y$ of $X^{(q)}$.  If each $q$-cell $\sigma$ of $Y$
has center point $p_\sigma$ then, by a uniform continuity argument, at
most finitely many of the components of $\Cont(f)$ meet the set
$f^{-1}[\{\,p_\sigma:\textrm{$\sigma$ a $q$-cell}\,\}]$.  For every
other component $U$, there is a homotopy rel $\bdy U$ from
$f\restriction\bar U$ to a map $g_U\colon\bar U\to X^{(q-1)}$.  Attach
these homotopies to the identity homotopy on the rest of $M$ to obtain
$H\colon f\homot g$.  This $H$ gives us the following:

\begin{lemma}\label{L:finiteComp}
  Any quasi-cellular $f\colon M\to X$ is nicely homotopic to $g\colon
  M\to X$ where $\Volcoh g=\Volcoh f$, $\VolAWP g=\VolAWP f$ when this
  expression makes sense, $\Voltr g\le\Voltr f$, and $\Cont(g)$ has
  finitely many components.
\end{lemma}

The statement for $\VolAWP$ is clear.  The $\Voltr$ statement follows
from the fact that $\Vol_\tau g=\Vol_\tau f$ for any $\tau$ where
$\Vol_\tau f$ is finite, which in turn follows from the $\VolAWP$
statement applied to each $q$-cell of $\tau$.  The $\Volcoh$ equivalence
follows from the commutative diagram
\begin{equation*}
  \xymatrix{
    \ckH^q(X,X\setminus e^\circ)
      \ar[r]^-{g^*} \ar@/_1pc/[rr]_{f^*=g^*} &
    \ckH^q(M,M\setminus g^{-1}[e^\circ])
      \ar[r]^{\iota^*} &
    \ckH^q(M,M\setminus f^{-1}[e^\circ])
  }
\end{equation*}
where $e$ is any $q$-cell of $X$ and $\iota$ is the canonical
inclusion (note that $\iota^*$ is a length-preserving embedding).

In fact, we can make lemma \ref{L:finiteComp} stronger.

\begin{lemma}\label{L:smoothComp}
  In lemma \ref{L:finiteComp} above, such a map $g$ exists where the
  boundary of $\Cont(g)$ is a smooth submanifold of $M$ and where, for
  each component $W$ of $\Cont(g)$, the restriction $g\restriction\bar
  W$ factors as $i\comp g_W$, where $i\colon D^q\to X$ canonically maps
  $D^q$ onto $\bar{\sigma(W)}$.
\end{lemma}

\begin{proof}
  Assume WLOG that $\Cont(f)$ has finitely many components.
  Consider each $q$-cell as a unit disk with the standard
  differential structure and geometry.
  By a mollification, $f$ is homotopic rel $f^{-1}[X^{(q-1)}]$ to a
  quasi-cellular map $g_1$ which is smooth on $\Cont(g_1)=\Cont(f)$.
  Clearly $\Volcoh g_1=\Volcoh f$, \emph{etc.}  For $x\in X^{(q)}$,
  let $\rho(x)=\dist(x,X^{(q-1)})\in[0,1]$.
  Then there exist
  arbitrarily small $\epsilon>0$ which are regular values of
  $\rho\comp g_1$.

  For a given $\epsilon$, compose $g_1$ with a homotopy of each
  $\sigma$ which deformation retracts an $\epsilon$-neighborhood of
  $\bdy\sigma$ onto $\bdy\sigma$ and scales the rest of $\sigma$
  uniformly, resulting in a map $g_2\colon M\to X$.  Then
  $\Cont(g_2)\subseteq\Cont(g_1)$.  If $\epsilon$ is chosen as above,
  then $\bdy\Cont(g_2)=(\rho\comp g_1)^{-1}[\epsilon]$ is a smooth
  submanifold of $M$.  For any component $W$ of $\Cont(g_2)$, the
  restriction $g_2\restriction\bar W$ factors through the map $i\colon
  D^q\to X$ as in the statement because $g_1\restriction\bar W$ maps
  to the slightly smaller closed ball of radius $1-\epsilon$.
  
  We do not yet know that
  $\Volcoh g_2\le\Volcoh g_1$, because we do not know that
  each component $V$ of $\Cont(g_1)$ contains at most one component
  of $\Cont(g_2)$.  We may achieve this by further homotopy.  First,
  by a uniform continuity argument, at most finitely many components
  of $V\cap\Cont(g_2)$ cover the center point of $\sigma(V)$; on the
  rest, we homotope $g_2$ to a map $g_3$ into $\bdy\sigma(V)$ as
  before.  Enumerate the components of $W=V\cap\Cont(g_3)$ as $W_1$,
  \dots, $W_n$, and suppose $n>1$.  We may homotope $g_3$ into a new
  map $h$ for which $V\cap\Cont(h)$ has $n-1$ components, iterate this
  construction to get at most one component, and repeat this iteration
  for each component $V$ to get a map $g$.  That $\Volcoh g=\Volcoh
  g_1$ may be seen by a commutative diagram similar to that in lemma
  \ref{L:finiteComp}.

  To construct $h$, let $T=(V\setminus W)^\circ$.  There is some
  component of $T$ which shares boundary points with both $W_1$ and
  $W_j$ for some $j>1$; otherwise the union of $\bar{W_1}$ and the
  adjacent components of $T$ is a nontrivial component of $V$, but $V$
  is by assumption connected.  Note that $\bdy W_1\cap\bdy
  W_j=\emptyset$.  Let $\gamma\colon[0,1]\to V$ be a path transverse
  to $\bdy W_1\cap\bdy T$ at 0 and transverse to $\bdy W_j\cap\bdy T$
  at 1, and where $\gamma[(0,1)]\subseteq T$.  Then $g_2$ is
  transformed into $h$ by a homotopy which moves into the cell
  $\sigma(V)$ on a neighborhood of $\gamma$ and is the identity away
  from $\gamma$.  The components of $V\cap\Cont(h)$ are those of
  $V\cap\Cont(g_2)$, except with $W_1$ and $W_j$ joined along
  $\gamma$, and by care in choosing the homotopy, the boundary of this
  joined component can be taken as smooth.
\end{proof}

Suppose $f$ is as in the conclusion of lemma~\ref{L:smoothComp}.  The
components of $\Cont(f)$ are manifolds with boundary, with the
boundaries disjoint; so we can express the volume as a sum
of degrees, as in our original motivation:
{\allowdisplaybreaks
\begin{align*}
  \Volcoh f &= \sum_\sigma \|f_\sigma^*([X^{(q)}])\| \\
  &= \sum_\sigma \bigl\| \bigl(f^*\colon \Zn\cong
  \ckH^q(X^{(q)},X^{(q)}\setminus\sigma^\circ)
  \to \ckH^q(M,M\setminus f^{-1}[\sigma^\circ])\bigr)[X^{(q)}] \bigr\| \\
  &= \sum_\sigma \sum_{j:f[W_j]\subseteq\sigma^\circ} \bigl\|f^*\colon
  \ckH^q(X^{(q)},X^{(q)}\setminus\sigma^\circ)\to
  \ckH^q(M,M\setminus W_j) \cong\Zn \bigr\| \\
  &= \sum_j \bigl\|f^*\colon
  \ckH^q(X^{(q)},X^{(q)}\setminus(\sigma(W_j))^\circ)\to
  \ckH^q(\overline{W}_j,\bdy W_j) \bigr\| \\
  &= \sum_j \bigl|\deg\bigl(f\colon (\overline{W}_j,\bdy W_j)\to
  (X^{(q)},X^{(q)}\setminus(\sigma(W_j))^\circ)\bigr)\bigr|.
\end{align*}}
The third line is by the direct sum decomposition
of $\ckH^q(M,M\setminus W)$ by the components of $W$.  The fourth
follows by excision and a homotopy equivalence.  The
last is easily seen through the universal coefficient theorem for
cohomology, which says in this case that the natural transformation
$\ckH^q\cong H^q\to \Hom(H_q(\text{---}),\Zn)$ is an isomorphism.

\begin{proof}[Proof of lemma \ref{L:AWPvscoh}]
  The second part of the lemma is easier.  The homotopy class of $f$
  is an element of $\pi_q(X^{(q)},X^{(q-1)},v)$ for some $v\in
  X^{(0)}$ (or is equivalent to an element of
  $\pi_1(X^{(1)}/X^{(0)},*)$).  This group has a standard generating
  set and a standard composition law; express $[f]$ in as few
  generators as possible, and let $f'$ be the map so obtained.  From
  inspection one sees that $\Volcoh f'=\VolAWP f'=\VolAWP f$.  There is
  a homotopy $H\colon f'\homot f$ where $H\restriction
  S^{q-1}\times[0,1]$ maps into $X^{(q-1)}$.  Now let $g\colon D^q\to
  X$ where $g(r,\theta)=f'(2r,\theta)$ for $r\le 1/2$ and
  $g(r,\theta)=H(\theta,2r-1)$ for $r\ge 1/2$.  $f'$ and $g$ have the
  same volumes, and $f$ is nicely homotopic to $g$.

  To prove the first statement, first take $q=1$, so that
  $f\colon(S^1,*)\to(X^{(1)},X^{(0)})$.  By lemma \ref{L:finiteComp},
  we may assume that $\Cont(f)$ has finitely many components.  By
  further homotopy we can eliminate any component which does not
  traverse an edge; this does not change $\VolAWP f$ or $\Volcoh f$.
  If any two adjacent components traverse the same edge in opposite
  directions, then yet further homotopy can eliminate these
  components, which lowers $\Volcoh f$ and keeps $\VolAWP f$ the same.
  After finitely many such reductions, we have a map which represents
  a reduced word in the generators of $K_1(X)$, and it is clear that
  $\VolAWP f = \Volcoh f$ is the number of remaining components.  This
  proves the case.

  Next, take $q\ge3$.  First adjust $f$ as in lemma~\ref{L:smoothComp}.
  Also assume $\pi_1(X)=0$; if this is not the case, lift $f$ to a map
  $\tilde f\colon M\to\tilde X$, which has the same volumes. Let $W=\Cont(f)$,
  and let $W_1$, \dots~$W_k$ be the components of $W$. Consider
  the commutative diagram:
  \begin{equation*}
    \xymatrix{
      \pi_q(X^{(q)},X^{(q-1)}) \ar[r]^-\sim &
      H_q(X^{(q)},X^{(q-1)}) \ar@{}[r]|(.7)\cong & \bigoplus_\sigma\Zn \\
      & H_q(D^q,D^q\setminus W) \ar@{}[r]|(.7)\cong \ar[u]^{f_*} &\Zn^k \\
      \pi_q(D^q,S^{q-1}) \ar[uu]^{f_*} \ar[r]^-\sim &
      H_q(D^q,S^{q-1}) \ar@{}[r]|(.7)\cong \ar[u] & \Zn
    }
  \end{equation*}
  $S^{q-1}$ and $X^{(q-1)}$ are simply connected, and the pairs
  $(D^q,S^{q-1})$ and $(X^{(q)},X^{(q-1)})$ are both
  $(q-1)$-connected, so by the Hurewicz isomorphism theorem (see
  \cite[\S 7.5]{ehS66}) the horizontal arrows are isomorphisms.  The
  isomorphism $H_q(D^q,D^q\setminus W)\cong\Zn^k$ may be seen as the
  product of the $H_q(D^q,D^q\setminus W_j)$ or as a direct sum
  diagram from the $H_q(\overline{W}_j,\bdy W_j)$.  The arrow from
  $H_q(D^q,S^{q-1})$ to $H_q(D^q,D^q\setminus W)$ is the identity on
  every factor of $\Zn^k$.  Thus
  {\allowdisplaybreaks
  \begin{align*}
    \VolAWP f &= |f_*[D^q]| \\
    &= \biggl| \sum_\sigma \bigl( f_*\colon H_q(D^q,S^{(q-1)})\to
    H_q(X^{(q)},X^{(q)}\setminus\sigma^\circ) \bigr)[D^q]
    \biggr| \\
    &\le \sum_\sigma \left|\bigl( f_*\colon H_q(D^q,S^{q-1})\to
    H_q(X^{(q)},X^{(q)}\setminus\sigma^\circ)
    \bigr)[D^q]\right| \\
    &= \sum_j \bigl\| f_*\colon H_q(\overline{W}_j,\bdy W_j)\to
    H_q(X^{(q)},X^{(q)}\setminus(\sigma(W_j))^\circ) \bigr\| \\
    &= \sum_j \bigl|\deg\bigl(f\colon (\overline{W}_j,\bdy W_j)\to
    (X^{(q)},X^{(q)}\setminus(\sigma(W_j))^\circ)\bigr)\bigr| \\
    &= \Volcoh f.
  \end{align*} }

  Finally, let $q=2$.  First adjust $f$ as in
  lemma~\ref{L:smoothComp}.  Next we show that, up to nice homotopy,
  we can assume that every component $W$ of $\Cont(f)$ is a disk.
  Initially, the boundary of any given $W$ has finitely many
  components, each of which is a smooth closed curve.
  $W$ is connected and bounded, so there must be a
  single component $\gamma\subseteq\bdy W$ which surrounds $W$ and all
  the other components $\gamma_i$ of $\bdy W$, and no $\gamma_i$ lies inside
  any other $\gamma_j$.

  Proceed as in the proof of lemma~\ref{L:smoothComp}, except ``in
  reverse'': connect each $\gamma_i$ to $\gamma$ by some path
  $\delta_i$ in $W$, then homotope $f$ to a new function $f'$ which is
  mostly the same as $f$, but maps into $X^{(q-1)}$ in a narrow strip
  around each $\delta_i$.  Then $\Cont(f')\subseteq\Cont(f)$ and each
  component of $\Cont(f)$ contains a unique component of $\Cont(f')$,
  which by earlier reasoning implies $\Volcoh f'=\Volcoh f$.
  Moreover, each component of $\Cont(f')$ has connected boundary,
  hence is a disk.

  Now choose a point $x\in\bdy D^2$ and draw $k-1$ loops starting at
  $x$, separating $D^2$ into $k$ regions, each containing a single
  component of $\Cont(f')$.  Each of these regions is topologically a
  disk, from which it becomes clear that
  \begin{equation*}
    [f'] = \gamma_1*(n_1\phi_{\sigma(1)}) + \dots +
           \gamma_k*(n_k\phi_{\sigma(k)})
  \end{equation*}
  where $\gamma_i$ is the composition of $f'$ with a path from $x$ to
  a point on the $i$th disk, $\sigma(i)$ is the cell to which this
  disk maps, and $n_i$ is the degree of the map between disks.  From
  there we see that
  \begin{equation*}
    \Volcoh f'=\sum_{i=1}^k |n_i| \ge \|[f']\| = \VolAWP f' = \VolAWP f.
    \qedhere
  \end{equation*}
\end{proof}

\begin{proof}[Proof of lemma \ref{L:trvscoh}]
  Again the second statement is easier.
  If $\Voltr f < \infty$, choose $\tau$ where $\Vol_\tau f = \Voltr f$. Apply
  lemma \ref{L:AWPvscoh} on each $q$-cell of $\tau$ to obtain $h$,
  nicely homotopic to $f$, where $\Volcoh h \le \Vol_\tau f = \Voltr f$.
  If $\Voltr f = \infty$, let $h=f$.
  
  For the first statement, first adjust $f$ as in lemma~\ref{L:smoothComp}.
  Then
  \begin{equation*}\Volcoh f = \sum_W |\deg f_W|\end{equation*}
  where $W$ ranges over
  the components of $\Cont(f)$ and $f_W\colon (\bar W,\bdy W)\to(D^q,S^{q-1})$
  as in the conclusion of lemma~\ref{L:smoothComp}. Let $\tau_W$
  triangulate $\bar W$ for each $W$, and let
  $\tau$ extend the $\tau_W$ to triangulate $M$.
  
  First we note that each $f_W$ is homotopic relative to $\bdy W$
  to a map $g_W$ which maps all of $M$, except for the interior of some
  $q$-cell $\Delta$ of $\tau_W$, to $S^{q-1}$. First consider the special case
  where $\tau_W$ consists of 2 $q$-cells $\Delta_1$ and $\Delta_2$ which
  intersect along a $(q-1)$-face $F$. It is not hard to see that $1_W\colon W\to W$
  is homotopic (relative to $\bdy W$) to a map which maps $\Delta_2$ onto $W$
  and $\Delta_1$ onto $\bdy \Delta_1\setminus F^\circ$. Composing this homotopy
  with $f_W$ gives $g_W$.
  
  For the general case, let $G$ be the undirected graph whose vertices are
  the $q$-cells of $\tau_W$ and where $\{\Delta,\Delta'\}$ is an edge iff
  $\Delta\cap\Delta'$ is an $(q-1)$-cell. Since $W$ is a connected manifold,
  $G$ is a connected graph. Let $T$ be a spanning tree for $G$. Proceed as
  follows: Let $\Delta_1$ be a leaf of $T$, and let $\Delta_2$ be the unique
  vertex where $\{\Delta_1,\Delta_2\}$ is an edge of $T$. Apply the above
  homotopy to $\Delta_1\cup\Delta_2$, then remove $\Delta_1$ and its edge
  from $T$. Repeat until $T$ has a single vertex. We have built a
  homotopy from $f_W$ to some $g_W$, and the key $q$-cell $\Delta$ is precisely
  the vertex remaining in $T$. By further homotopy, we can make $g_W$
  smooth on $g^{-1}[(D^q)^\circ]$ and cover each point of $(D^q)^\circ$ minimally.
  
  We now build the homotopy from $f$ to $g$. On each $W$, compose the homotopy
  from $f_W$ to $g_W$ with the characteristic map of the disk into $X$ to obtain
  a homotopy from $f\restriction W$ to $g\restriction W$. Outside of $\Cont(f)$,
  take the constant homotopy. Then
  \begin{equation*}
    \Voltr g \le \Vol_\tau g = \sum_W |\deg f_W| = \Volcoh f
  \end{equation*}
  and the rest of the claims hold as well.
\end{proof}

\section{Filling volumes and profiles}\label{S:filling}

Recall that $\mathbf{I}_t(X)$ is the space of integral $t$-dimensional
currents in $X$, while $C^{0,1}(M,X)$ is the space of Lipschitz maps
from $M$ to $X$.

\begin{definition}\label{D:fillvol}
  Let $X$ be a Riemannian manifold.  If $S\in\mathbf{I}_{q-1}(X)$,
  then the \emph{current filling volume} of $S$ is
  \begin{equation*}
    \FVcur(S) := \inf\,\{\,\mathbf{M}(T) : T\in\mathbf{I}_q(X),\,\bdy T=S\,\}.
  \end{equation*}
  If $f\in C^{0,1}(\bdy M,X)$, then the \emph{Lipschitz filling
    volume} of $f$ through $M$ is
  \begin{equation*}
    \FVmet^{X,M}(f) := \inf\,\{\,\Volmet(h) : h\in C^{0,1}(M,X),\,\bdy h=f\,\}.
  \end{equation*}
  Now let $X$ be a CW complex.  If $c$ is a $(q-1)$-chain in $X$, then
  the \emph{chain filling volume} of $c$ is
  \begin{equation*}
    \FVch(c) := \inf\,\{\,\|b\| : b\in C_q(X),\,\bdy b=c\,\}.
  \end{equation*}
  If $f\colon \bdy M\to X$ is quasi-cellular, then the
  \emph{cellular filling volume} of $f$ through $M$ is
  \begin{align*}
    \FVcell^{X,M}(f) :=& \inf\,\{\,\Volcoh(g) : g\colon M\to X\textrm{ quasi-cellular, }\bdy g=f\,\} \\
    =& \inf\,\{\,\Voltr(g) : g\colon M\to X\textrm{ quasi-cellular, }\bdy g=f\,\}\\
    =& \inf\,\{\,\Vol(g) : g\colon M\to X\textrm{ admissible, }\bdy g=f\,\}; \\
    \intertext{in the case $M=D^q$, also}
    =& \inf\,\{\,\VolAWP(g) : g\colon D^q\to X\textrm{ quasi-cellular, }\bdy g=f\,\}.
  \end{align*}
\end{definition}
The various definitions for $\FVcell$ are equivalent by an application
of lemmas~\ref{L:AWPvscoh} and~\ref{L:trvscoh}. Note that $\FVcell^{X,D^q}$
is the filling volume defined on \cite[p.~87]{jmA99}.

\begin{lemma}\label{L:niceHomotopy}
  If $f,\,g\colon \bdy M\to X$ are quasi-cellular and nicely homotopic, then
  \begin{equation*}\FVcell^{X,M}(f)=\FVcell^{X,M}(g).\end{equation*}
\end{lemma}

\begin{proof}
  Let $F\colon M\to X$ be quasi-cellular with $\bdy F=f$, and let $H\colon f\homot g$
  be a nice homotopy. There is a homeomorphism
  \begin{equation*}
  \phi\colon M \cong M' = 
      \bigl[M \sqcup (\bdy M\times[0,1])\bigr]\big/\{\, x \sim (x,0) : x\in\bdy M\,\}
  \end{equation*}
  because $\bdy M$ has a collar neighborhood; $F$ and $H$ together form a
  continuous map $F'$ on $M'$, and $G=F'\comp\phi$
  is a quasi-cellular map filling $g$ with $\Volcoh G = \Volcoh F$. This proves
  $\FVcell^{X,M}(f) \ge \FVcell^{X,M}(g)$; the reverse is similar.
\end{proof}

\begin{definition}
  Let $X$ be a Riemannian manifold. The \emph{current profile} of $X$ in dimension
  $q$ is the function $\Phicur^{X,q}\colon [0,\infty) \to [0,\infty]$ where
  \begin{equation*}
    \Phicur^{X,q}(v) := \sup\,\{\,\FVcur(\bdy T) : T\in\mathbf{I}_q(X),\,\mathbf{M}(\bdy T)\le v\,\}.
  \end{equation*}
  The \emph{metric profile} of $X$ for $M$ is the function
  $\Phimet^{X,M}\colon[0,\infty)\to[0,\infty]$ where
  \begin{equation*}
      \Phimet^{X,M}(v) := \sup\,\{\,\FVmet^{X,M}(\bdy h) : h\in C^{0,1}(M,X),\,\Volmet(\bdy h)\le v\,\}.
  \end{equation*}
  Now let $X$ be a CW complex. The \emph{chain profile} of $X$ in dimension
  $q$ is the function $\Phich^{X,q}\colon\Nn\to \Nn\cup\{\infty\}$ where
  \begin{equation*}
    \Phich^{X,q}(n) := \sup\,\{\,\FVch(\bdy b) : b\in C_q(X),\,\|\bdy b\|\le n\,\}.
  \end{equation*}
  The \emph{cellular profile} of $X$ for $M$ is the function
  $\Phicell^{X,M}\colon\Nn\to\Nn\cup\{\infty\}$ where
  \begin{align*}
    \Phicell^{X,M}(n) :=& \sup\,\{\,\FVcell^{X,M}(\bdy f) : f\colon M\to X\text{ quasi-cellular, }
      \Volcoh \bdy f\le n\,\} \\
    =& \sup\,\{\,\FVcell^{X,M}(\bdy f) : f\colon M\to X\text{ quasi-cellular, }
      \Voltr \bdy f\le n\,\}\\
    =& \sup\,\{\,\FVcell^{X,M}(\bdy f) : f\colon M\to X\text{ quasi-cellular, }
      \Vol \bdy f\le n\,\}; \\
    \intertext{in the case $\bdy M= S^{q-1}$, also}
     =& \sup\,\{\,\FVcell^{X,M}(\bdy f) : f\colon M\to X\text{ quasi-cellular, }
      \VolAWP \bdy f\le n\,\}.
  \end{align*}
\end{definition}
The definitions for $\Phicell^{X,M}$ are equivalent by an application of
lemmas~\ref{L:AWPvscoh}, \ref{L:trvscoh}, and~\ref{L:niceHomotopy}.
Note that $\Phicell^{X,D^q}$ is the $q$-dimensional Dehn function
defined on \cite[p.~90]{jmA99}.

These functions are the generalized Dehn functions we wish to study. Our positive
results will state that two such functions are equal, or equal almost everywhere,
or equivalent as follows.

\begin{definition}
  Let $f,\,g\colon[0,\infty)\to[0,\infty]$ be weakly increasing. $f$ is
  \emph{quasi-bounded} by $g$, written $f\qbound g$, if there exist $A$,~$B>0$
  and $C$,~$D\ge 0$ where
  \begin{equation*}
    f(x) \le A\cdot g(Bx) + Cx + D\qquad\text{for all $x\ge 0$.}
  \end{equation*}
  We say $f$ and $g$ are \emph{quasi-equivalent,} written $f\qequiv g$, if
  $f\qbound g\qbound f$.
\end{definition}
A routine check shows that $\qbound$ is a preorder, making $\qequiv$ an
equivalence relation. We extend these concepts to functions $f$ on $\Nn$ by
extending $f$ to $[0,\infty)$ so that $f(x) = f(\lfloor x\rfloor)$. Ideas similar to $\qbound$
and $\qequiv$ appear throughout the literature.

\section{Metric and cellular profiles}\label{S:metvscell}

For this section, let $\pi\colon X\to Y$ be a covering map, where $Y$ is a compact
Lipschitz neighborhood retract (CLNR; say a Riemannian manifold)
with a Lipschitz triangulation $\tau$. Both the metric structure and the triangulation
may be lifted through $\pi$ to $X$, so that all of the profiles in
section~\ref{S:filling} are defined for $X$. As we might expect, and will
prove in this section, the profiles corresponding to a manifold $M$ are equivalent,
as are the homological profiles in each dimension $q$. Briefly, this occurs because
the two notions of volume for a quasi-cellular map are roughly equivalent, and
any Lipschitz map $f$ may be deformed into a quasi-cellular map without changing
$\Volmet f$ by more than a constant factor.

\begin{theorem}\label{T:metequivcell}
  $\Phicur^{X,q}\qequiv\Phich^{X,q}$ for any $q\ge 2$ and
  $\Phimet^{X,M}\qequiv\Phicell^{X,M}$ for any $M$.
\end{theorem}

Thus there is a ``metric independence'' result:
\begin{corollary}
  If $Y$ is a closed Riemannian manifold [with boundary], then the profiles
  $\Phimet^{X,M}$ and $\Phicur^{X,q}$ are independent of the metric on $Y$,
  up to quasi-equivalence.
\end{corollary}

We spend the rest of the section proving theorem~\ref{T:metequivcell}. First note
that any triangular $q$-chain in $X$ is also a $q$-current, and there are constants
$0<C\le D$ (the minimum and maximum volume of a $q$-simplex,
respectively) where
\begin{equation*}
C\Volch T\le \mathbf{M}(T) \le D\Volch T
\end{equation*}
for all triangular $q$-chains $T$.  A similar statement for
 quasi-cellular maps from $M$ is harder to state and prove.

\begin{lemma}\label{L:gEqCoh}
  For every quasi-cellular map
  $f\colon M\to X$, $C\Volcoh f\le\Volmet f$ if $f$ is Lipschitz.
  Moreover, there is a Lipschitz map $g\simeq f$ rel $\bdy M$ so that
  $\Volcoh g=\Volcoh f$ and $\Volmet g\le D\Volcoh g$.
\end{lemma}

\begin{proof}
  First assume for simplicity that $f$ is $C^1$ on $\Cont(f)$.
  We start with the area formula
  \begin{equation*}
    \Volmet f = \sum_{\sigma^q} \int_\sigma N(f,y)\,dh(y).
  \end{equation*}
  If $M$ is a measurable subset of $\Rn^q$, this follows from case (2)
  of \cite[theorem 3.2.3]{hF69}; for general $M$, use a
  partition-of-unity argument.  Now almost everywhere in any cell
  $\sigma$, $y$ is a regular value of $f$ by Sard's theorem.  If
  $f(x)=y$ for $y$ regular, then $f$ has a local degree $\deg_x f=\pm
  1$ at $x$, and for each component $W$ of $f^{-1}[\sigma]$
  \begin{equation*}
    |\deg (f\restriction\bar W)|
    =\biggl|\sum_{x\in W,f(x)=y} \deg_x f \biggr|
    \le N(f\restriction W,y).
  \end{equation*}
  After summing over components $W$, we have 
  \begin{equation*}
    \|f_{(\sigma)}^*([M])\|\le N(f,y).
  \end{equation*}
  Integrate over the $q$-cell $\sigma$ to obtain
  \begin{equation*}
    C\|f_{(\sigma)}^*([M])\|\le \int_\sigma N(f,y)\,dh(y)
  \end{equation*}
  and sum over all cells $\sigma$ to obtain $C\Volcoh f \le \Volmet f$.

  For general $f\in C^{0,1}(M,X^{(q)})$, use \cite[theorem
  3.1.15]{hF69}, which implies that there is a $C^1$ function $g$ with
  the same Lipschitz constant as $f$ so that $g=f$ except on a set $S$
  where $mS$ is arbitrarily small (assume a background metric).  If we
  only modify $f$ at points $x$ where $\dist(f(x),X^{(q-1)})>\epsilon>0$,
  and make $mS$ small enough, we can guarantee
  that $\Volcoh f$ does not change, and obtain a lower bound $\Volmet f\ge
  C(1-\epsilon)^q\Volcoh f-\|f\|_{0,1}^q mS$.  Now take $\epsilon$,~$mS\to0$.

  For the second part, choose $g\homot f$ and a suitable triangulation of
  $M$ as in lemma~\ref{L:trvscoh}.  We
  need only check the last condition.  Restrict $g$ to a single
  simplex $\Delta$.  If $g[\Delta]\subseteq X^{(q-1)}$, then
  $\Volmet(g\restriction\Delta)=0$.  Otherwise $g[\Delta]$ meets the
  interior of at most one cell $\sigma$, and we have
  \begin{equation*}
    \Volmet(g\restriction\Delta) = \int_{\sigma^\circ} N(g,y)\,dh(y)
    = |g_*([\Delta])|\Vol \sigma \le D|g_*([\Delta]).
  \end{equation*}
  Sum over the simplices $\Delta$ to obtain $\Volmet g\le D\Volcoh g$.
\end{proof}

Next we must show that any map can be made quasi-cellular without
changing its volume by more than a constant factor. In other words,
we need the Deformation Theorem, except with functions in place of
currents.

\begin{lemma}\label{L:mapDeformation}
  There exists a constant $C$ where, for all $f\in C^{0,1}(M,X)$,
  there exists a Lipschitz homotopy $H\colon f\homot f'$ where
  $f'[M]\subseteq X^{(q)}$ and where
  \begin{align*}
    \Volmet f' &\le C\Volmet f,\\
    \Volmet \bdy f' &\le C\Volmet \bdy f,\\
    \text{and }\Volmet H &\le C\Volmet f.
  \end{align*}
  Moreover, whenever $f(x)\in X^{(q)}$, we have $H(t,x)=f(x)$,
  and in particular $f'(x)=f(x)$.
\end{lemma}

\begin{proof} Adapted from \cite{dbaE92}, starting on page 223, which is
  itself adapted from the classic proofs in \cite{hF69} and \cite{lS84}.
  WLOG assume the simplices of $X$ are standard.
  Suppose $f[M]\subseteq X^{(k)}$ with $k>q$.  It suffices to show that
  we can deform $f$ to
  $f'$ where $f'[M]\subseteq X^{(k-1)}$, the estimates above hold for
  some constant $C$, and those points already lying in $X^{(k-1)}$ are
  unmoved; from there we simply iterate the deformation to
  $X^{(k-2)}$, \emph{etc.}, until our map lies in $X^{(q)}$.  Further,
  we can perform this deformation separately on the interior of each
  $k$-simplex $\Delta$ and glue the results together by local
  finiteness of $(X,\tau)$.

  Thus, fix $\Delta$ and let $S=f^{-1}[\Delta^\circ]$.  Let $u_0$ be the barycenter of
  $\Delta$, and let $r=\dist(u_0,\bdy\Delta)$. For every $u\in B:=
  B(u_0,r/4)$, let $\pi_u\colon \Delta\setminus\{u\}\to\Delta\setminus
  B(u,r/2)$ be the identity outside the ball $B(u,r/2)$ and radial
  projection inside.  The function $f_u = \pi_u\comp(f\restriction S)$
  is defined for almost all $u$, and
  \begin{align*}
    \Volmet f_u &= \int_{S\setminus f^{-1}[B(u,r/2)]} f_u(x)\,dx
      +  \int_{f^{-1}[B(u,r/2)]} f_u(x)\,dx \\
    &\le \Volmet(f\restriction S) + \int_{f^{-1}B(u,r/2)}
      \frac{J\!f(x)(r/2)^q}{|f(x)-u|^q}\,dx.
  \end{align*}
  If we integrate over all $u\in B$, we have
  {\allowdisplaybreaks
  \begin{align*}
    \int_B (&\Volmet f_u)\,du \\
    &\le \Vol(B)\Volmet(f\restriction S) +
    \int_B\int_{f^{-1}B(u,r/2)}
    \frac{J\!f(x)(r/2)^q}{|f(x)-u|^q}\,dx\,du \\
    &\le c_1\Volmet(f\restriction S) \\ & \quad +
    \left(\frac{r}{2}\right)^q \int_{f^{-1}B(u_0,3r/4)}
    dx\left[J\!f(x)
      \int_{B(f(x),r/2)} |u-f(x)|^{-q}\chi_B(u)\,du\right] \\
    &\le c_1\Volmet(f\restriction S) + \left(\frac{r}{2}\right)^q
    \left(\int_S J\!f(x)\,dx\right)
    \left(\int_{B(0,r/2)} |v|^{-q}\,dv\right) \\
    &\le c\Volmet(f\restriction S)
  \end{align*} }
  since $\int |v|^{-q}\,dv$ is finite over balls in $\Rn^k$, $k>q$.
  (In this case $c$ is dependent only on $k$ and $q$.)  It follows
  that, if we define for $v>0$
  \begin{equation*}
    A_v := \{\,u\in B:\Volmet f_u > v\Volmet(f\restriction S)\,\},
  \end{equation*}
  then $mA_v \le c/v$.  We can perform similar estimates for $\bdy
  f$ and for the radial homotopies $h_u$; if
  \begin{align*}
    B_v &:= \{\,u\in B:\Volmet (\bdy f)_u > v\Volmet(\bdy f\restriction
    (S\cap\bdy M))\,\} \\
    \intertext{and}
    C_v &:= \{\,u\in B:\Volmet h_u > v\Volmet(f\restriction S)\,\},
  \end{align*}
  then $mB_v \le c'/v$ and $mC_v\le c''/v$ for some $c'$, $c''$
  depending on $k$ and $q$.

  Thus, if we fix $v>(c+c'+c'')/\Vol(B)$ (once again $v=v(k,q)$), then
  there is some $u\in B$ where $\Volmet f_u\le v\Volmet(f\restriction S)$,
  $\Volmet (\bdy f)_u \le v\Volmet(\bdy f\restriction S)$, and $\Volmet h_u
  \le v\Volmet(f\restriction S)$.  Since $\dist(x_0,f_u[S])\ge r/4$, we
  may further deform $f_u$ to $f_u'\colon S\to\bdy\Delta$ by a radial
  homotopy, and all the volumes are appropriately bounded.
\end{proof}

\begin{lemma}\label{L:mapDeformation2}
  There is a constant $C$ where, for $f\in C^{0,1}(M,X)$, there is a
  Lipschitz homotopy $H\colon f\homot f'$ where $f'$ is quasi-cellular and 
  \begin{align*}
    \Volmet f' &\le C(\Volmet f+\Volmet \bdy f),\\
    \Volmet \bdy f' &\le C\Volmet \bdy f,\\
    \text{and }\Volmet H &\le C\Volmet f.
  \end{align*}
\end{lemma}

\begin{proof}
  First produce a homotopy $H'\colon f\homot f''$ as in lemma
  \ref{L:mapDeformation}, then produce a homotopy $J\colon\bdy f''\to
  j$, also as in lemma \ref{L:mapDeformation}.  We obtain $f'$ by
  attaching $J$ to $f''$ by a collar neighborhood, which yields the
  estimate on $\Volmet f'$.  $\bdy f'=j$ gives the estimate on $\Volmet
  f'$, and $f''\homot f'$ by a homotopy which lives entirely inside
  $f'[M]$ and therefore has zero volume; so $\Volmet H=\Volmet H'$ gives
  the last estimate.
\end{proof}

We can probably strengthen lemma~\ref{L:mapDeformation2} to remove the
$\Volmet\bdy f$ in the first estimate, \emph{a la} \cite{lS84}.

\begin{lemma}\label{L:FVgEqCoh}
  There exist constants $0<C\le D$ where, for all $f\colon\bdy M\to X^{(q-1)}$
  quasi-cellular and Lipschitz,
  \begin{equation*}
    C\FVcell^{X,M}(f) \le \FVmet^{X,M}(f) \le D\FVcell^{X,M}(f).
  \end{equation*}
\end{lemma}

\begin{proof}
  For the first inequality, let $h\in C^{0,1}(M,X)$ where $\bdy h=f$.
  By lemma \ref{L:mapDeformation}, there is a map $h'\in
  C^{0,1}(M,X^{(q)})$ where $\bdy h'=f$ and $\Volmet h'\le C_1\Volmet h$.
  By lemma \ref{L:gEqCoh}, there is a constant $c$ where $c\Volcoh
  h'\le \Volmet h'$.  Thus let $C=c/C_1$.  The second inequality follows
  directly from lemma \ref{L:gEqCoh}.
\end{proof}

The corresponding results for currents carry over directly. We state them
without proof.

\begin{lemma}\label{L:currentDeformation}
  There is a constant $C$ where, for all $T\in\mathbf{I}_q(X)$,
  there exist currents $P$, $S\in\mathbf{I}_q(X)$ and
  $R\in\mathbf{I}_{q+1}(X)$ where $T-P=\bdy R+S$, $\supp P\subseteq
  X^{(q)}$, and
  \begin{align*}
    \mathbf{M}(P) &\le C\mathbf{M}(T),&
    \mathbf{M}(\bdy P) &\le C\mathbf{M}(\bdy T),\\
    \mathbf{M}(R) &\le C\mathbf{M}(T),&
    \mathbf{M}(S) &\le C\mathbf{M}(\bdy T).
  \end{align*}
  If $\supp \bdy T \subseteq X^{(q)}$, then $\bdy S=0$ (\emph{i.e.},
  $\bdy P=\bdy T$).
\end{lemma}

\begin{lemma}\label{L:currentDeformation2}
  As above, but the conclusion includes $\supp\bdy P\subseteq
  X^{(q-1)}$ and the first volume estimate becomes
  \begin{equation*}
    \mathbf{M}(P) \le C(\mathbf{M}(T) + \mathbf{M}(\bdy T)).
  \end{equation*}
\end{lemma}

\begin{lemma}\label{L:FVcurEqCh}
  There exist constants $0<C\le D$ where, for any $q$-cycle $T$
  of $X$,
  \begin{equation*}
    C\FVch(T) \le \FVcur(T) \le D\FVch(T).
  \end{equation*}
\end{lemma}

\begin{proof}[Proof of theorem~\ref{T:metequivcell}]
  We prove the second statement.
  Let $f\colon\bdy M\to X$ be quasi-cellular; we may perturb it slightly
  to make it Lipschitz without changing $\Volcoh f$. Then by lemmas
  \ref{L:gEqCoh} and \ref{L:FVgEqCoh}
  \begin{equation*}
    \FVcell^{X,M}(f) \le c\FVmet^{X,M}(f)\le c\Phimet^{X,M}(\Volmet f)
    \le c\Phimet^{X,M}(D\Volcoh f).
  \end{equation*}
  It follows that for all $n$,
  \begin{equation*}
    \Phicell^{X,M}(n) \le c\Phimet^{X,M}(Dn).
  \end{equation*}
  Now let $f\in C^{0,1}(\bdy M,X)$.  By lemma \ref{L:mapDeformation}
  there is a Lipschitz homotopy $H\colon f\homot f'$ where $f'$ is
  quasi-cellular, $\Volmet f'\le C\Volmet f$, and $\Volmet H\le C\Volmet f$.
  We can attach $H$ to any filling function for $f'$ to get a filling
  function for $f$, so that
  \begin{align*}
    \FVmet^{X,M}(f) &\le \FVmet^{X,M}(f') + C\Volmet f
    \le D\FVcell^{X,M}(f') + C\Volmet f \\
    &\le D\Phicell^{X,M}(\Volcoh f') + C\Volmet f
    \le D\Phicell^{X,M}(c\Volmet f') + C\Volmet f \\
    &\le D\Phicell^{X,M}(cC\Volmet f) + C\Volmet f
  \end{align*}
  and for any $x\ge0$,
  \begin{equation*}
    \Phimet^{X,M}(x) \le D\Phicell^{X,M}(cCx) + Cx.
  \end{equation*}
  The first statement is proved similarly.
\end{proof}

\section{Profiles and $q$-homotopy equivalence}\label{S:homotopy}

In the usual proofs that the Dehn function $\delta^k$ of a group $G$ is well-defined
up to equivalence, one starts with two $K(G,1)$s $X$ and $Y$ with finite $(k+1)$-skeleta.
Without loss of generality, $Y$ is a subcomplex of $X$, in fact a deformation retraction.
By choosing the homotopy (rel $X$) between $1_X$ and $r\colon X\to Y$ to be cellular,
one obtains a way to deform maps $M\to\tilde X$ into maps $M\to\tilde Y$ while changing their
volume by at most a constant factor. This allows one to show that
$\delta_X^k\qequiv\delta_Y^k$.

This argument works just as well for maps from any $M^q$ into
$\tilde X$ and $\tilde Y$. Also,
$X$ and $Y$ may be any two homotopy equivalent complexes
with finite $q$-skeleta, not necessarily $K(G,1)$s; and we need not
use the universal covers, only covers corresponding to the same
subgroup of $\pi_1(X)\cong\pi_1(Y)$.
Add to this that the cells in dimensions higher
than $q$ are irrelevant, and we have the following.

\begin{theorem}\label{T:qConnected}
Let $q\ge 2$, let $X$ and $Y$ be connected CW complexes with $X^{(q)}$ and
$Y^{(q)}$ finite, and suppose there is a continuous map $f\colon Y\to X$ where
the induced map $f_*\colon \pi_t(Y,*)\to\pi_t(X,f(*))$ is an isomorphism for $1\le t<q$. Let
$\tilde Y$,~$\tilde X$ be covering spaces of $Y$ and $X$ corresponding to subgroups
$G\subseteq\pi_1(Y,*)$ and $f_*[G]\subseteq\pi_1(X,f(*))$ respectively.
Then $\Phicell^{\tilde X,M}\qequiv\Phicell^{\tilde Y,M}$ for all $q$-dimensional $M$, and
$\Phich^{\tilde X,q}\qequiv\Phich^{\tilde Y,q}$.
\end{theorem}

The special case where $M=D^q$ (which makes the covering spaces irrelevant)
and $\pi_t(X)=\pi_t(Y)=0$ for $1<t<q$ is an easy corollary to the major theorems
in~\cite{jmA99}.

To see that isomorphic fundamental groups are not sufficient (as
they are in~\cite{jmA99}), consider for each $n\ge 1$
the spaces $X=S^{2n}$ (with one 0-cell and one $2n$-cell) and
$Y=S^{2n}\sqcup_\alpha D^{4n}$, where $\alpha\colon S^{4n-1}\to S^{2n}$
represents an element $a\in\pi_{4n-1}(S^{2n})$ with infinite order. Then
$\Phicell^{X,D^{4n}}\equiv0$, since any map $D^{4n}\to S^{2n}$ has zero
volume. By contrast, $\Phicell^{Y,D^{4n}}\equiv\infty$, as any map
$S^{4n-1}\to Y^{(2n)}=S^{2n}$ representing $ka$ has volume 0 and filling
volume $|k|$.

\begin{proof}[Proof of theorem~\ref{T:qConnected}]
We show that $\Phicell^{\tilde X,M}\qequiv\Phicell^{\tilde Y,M}$; the reasoning for
$\Phich$ is analogous.
WLOG assume $f$ is cellular, and an inclusion (consider the mapping cylinder~$M_f$).
By adjoining $(q+1)$-cells to $X$, we may kill $\pi_q(X)$ without changing
$\Phicell^{\tilde X,M}$ or $\Phich^{\tilde X,q}$, so also assume $f$ is $q$-connected.
In other words, we reduce to the case of a $q$-connected pair $(X,A)$ with $X^{(q)}$ finite,
$\pi\colon \tilde X\to X$ a covering space, and $\tilde A=\pi^{-1}[A]$ the corresponding
covering space of $A$.

Let $j\colon (X^{(q)},A)\to (X,A)$ be the inclusion map.  $(X,A)$ is
$q$-connected, so there is a (cellular) homotopy $h\colon j\homot g$
to some cellular map $g\colon X^{(q)}\to A$, and both $g$ and $h$
fix $A$.
Assume WLOG that $h\restriction\sigma$ is admissible
for each $q$-cell $\sigma$ of $X^{(q)}\times I$.
Let $K := \sup_\sigma \Vol(h\restriction\sigma)$,
which is
finite because only finitely many $\sigma$ exist. Assume WLOG
 that $g\restriction\sigma$ is
admissible for each $(q-1)$-cell $\sigma$ of $X$, and let
$L := \sup_\sigma \Vol(g\restriction\sigma)$,
which is also finite.
$g$ lifts to a map $\tilde g\colon\tilde X^{(q)}\to\tilde A$ which
fixes $\tilde A$,
and $h$ lifts to a homotopy $\tilde h\colon\tilde X^{(q)}\times I\to \tilde X$
from the inclusion $\tilde\jmath$ to $\tilde g$, also fixing $\tilde A$.

To establish $\Phicell^{\tilde X,M}\qbound \Phicell^{\tilde A,M}$, fix $n$ and let
$\omega=\bdy\psi$ where $\psi\colon M\to\tilde X$ and $\omega$ are both
admissible and
$\Vol \omega \le n$.
Then $\omega':=\tilde g\comp\omega\colon\bdy M\to\tilde A$ is 
admissible, $\omega'=\bdy(\tilde g\comp\psi)$ in $\tilde A$, and
$\Vol \omega'\le Ln$.
Similarly $\omega'':=\tilde h\comp(\omega\times 1_I)\colon\bdy M\times I\to X^{(q)}$
is an admissible homotopy from $\omega$ to $\omega'$, and
$\Vol \omega''\le Kn$. Let $\phi'\colon M\to\tilde A$ fill $\omega'$
where $\Vol\phi'=\FVcell^{\tilde A,M}(\omega')$. By attaching $\phi'$ to $\omega''$, we obtain
an admissible filling map for $\omega$, and we see that
\begin{equation*}
\FVcell^{\tilde X,M}(\omega) \le \Vol(\phi')+\Vol(\omega'')
\le\FVcell^{\tilde A,M}(\omega')+Kn\le\Phicell^{\tilde A,M}(Ln)+Kn.
\end{equation*}
We take the supremum over all such $\omega$ to see that
\begin{equation*}
\Phicell^{\tilde X,M}(n) \le \Phicell^{\tilde A,M}(Ln)+Kn,
\end{equation*}
so in particular $\Phicell^{\tilde X,M}\qbound\Phicell^{\tilde A,M}$.

To establish $\Phicell^{\tilde A,M}\qbound\Phicell^{\tilde X,M}$, fix $n$ and let
$\omega=\bdy\psi$ where $\psi\colon M\to\tilde A$ and $\omega$ are admissible and
$\Vol\omega\le n$. Let $\phi\colon M\to\tilde X$ fill $\omega$ where
$\Vol\phi=\FVcell^{\tilde X,M}(\omega)$. Then
$\phi'=g\comp\phi\colon M\to\tilde A$ is admissible and
$\FVcell^{\tilde A,M}(\omega)\le\Vol\phi'\le K\Vol\phi\le K\Phicell^{\tilde X,M}(n)$.
Taking the supremum over all such $\omega$, we see that
\begin{equation*}
\Phicell^{\tilde A,M}(n) \le K\Phicell^{\tilde X,M}(n),
\end{equation*}
and therefore $\Phicell^{\tilde A,M}\qbound\Phicell^{\tilde X,M}$.
\end{proof}

There is a metric version of theorem~\ref{T:qConnected}, which is an easy consequence
of it and of theorem~\ref{T:metequivcell}:

\begin{corollary}
Let $q\ge 2$, let $X$ and $Y$ be connected triangulable CLNRs (say connected
compact Riemannian manifolds), and suppose there
is a continuous map $f\colon Y\to X$ where $f_*\colon \pi_t(Y,*)\to\pi_t(X,f(*))$ is an
isomorphism for $1\le t<q$. Let $\tilde Y$,~$\tilde X$ be covering spaces of $Y$ and $X$
corresponding to subgroups $G\subseteq\pi_1(Y,*)$ and $f_*[G]\subseteq\pi_1(X,f(*))$
respectively. Then $\Phicur^{\tilde X,q}\qequiv\Phicur^{\tilde Y,q}$,
and $\Phimet^{\tilde X,M}\qequiv\Phimet^{\tilde Y,M}$ for all $q$-dimensional $M$.
\end{corollary}

\section{Further questions}\label{S:questions}

Intuitively, the large-scale geometry of a universal covering space $\tilde X$
with $X$ compact should be captured by $\pi_1(X)$. Thus there should be a
natural way of taming the infinities noted in section~\ref{S:homotopy} such that
the ``tamed'' versions of $\Phicell^{\tilde X,M}$ \emph{etc.} will depend only on $\pi_1(X)$
up to quasi-equivalence. We do not see any way to do so; on the other hand,
we also do not know of any compact spaces $X$, $Y$ where $\pi_1(X)\cong\pi_1(Y)$,
$\Phicell^{\tilde X,M}(n)<\infty$ and $\Phicell^{Y,M}(n)<\infty$ for all $n$, but
$\Phicell^{\tilde X,M}\not\qequiv\Phicell^{\tilde Y,M}$.

\bibliographystyle{hplain}
\bibliography{database}

\end{document}